\documentclass[11pt, final]{amsart}
\usepackage{amssymb}

\theoremstyle{plain}
\newtheorem{theorem}{Theorem}

\newtheorem{lemma}{Lemma}

\theoremstyle{remark}
\newtheorem*{remark}{Remark}

\def\E{{\mathcal E}(H)}
\def\la{\langle}
\def\N{\mathbb N}

\def\R{\mathbb R}
\def\ra{\rangle}
\def\rng{\operatorname{rng}}

\begin{document}
\title[]{Characterizations of the automorphisms of Hilbert
space effect algebras}
\author{LAJOS MOLN\'AR}
\address{Institute of Mathematics and Informatics\\
         University of Debrecen\\
         4010 Debrecen, P.O.Box 12, Hungary}
\email{molnarl@math.klte.hu}
\thanks{  This research was supported by the
          Hungarian National Foundation for Scientific Research
          (OTKA), Grant No. T030082, T031995, and by
          the Ministry of Education, Hungary, Reg.
          No. FKFP 0349/2000}
\date{\today}
\begin{abstract}
In this paper we characterize the automorphisms of Hilbert space
effect algebras by means of their preserving properties which
concern certain relations and quantities appearing in quantum
measurement theory.
\end{abstract}
\maketitle

\section{Introduction and statement of the results}

The concept of effects plays fundamental role in the mathematical
description of quantum measurement (for detailed explanations see
Introduction in \cite{BGL} or $\S$ 1 in \cite{Kraus}). In the
Hilbert space framework, the set $\E$ of all effects on a complex
Hilbert space $H$ is just the operator interval $[0,I]$ of all
positive operators on $H$ which are bounded by the identity $I$.
The set $\E$ can be equipped with several algebraic operations and
relations which all have physical content. Hence we have different
algebraic structures on the same set $\E$. The investigation of
the morphisms of these structures was initiated by Ludwig (for
explanation and results see Chapters V. and VI. in \cite{Ludwig}).
We recall one of his fundamental results in this direction as
follows. First, there is a natural partial order $\leq$ on $\E$
which is induced by the usual order between selfadjoint operators
on $H$. Next, there is a kind of orthocomplementation $\perp:
E\mapsto I-E$ on $\E$ (cf. \cite{BGL}, p. 25).  Now, Ludwig's
result \cite[Section V.5.]{Ludwig} (also see \cite{CVLL})
describes the ortho-order automorphisms of $\E$ (that is, the
automorphisms of $\E$ with respect to the relation $\leq$ and the
operation $\perp$) in the following way: If $\dim H\geq 3$, then
every ortho-order automorphism $\phi$ of $\E$ is of the form
\[
\phi(E)=UEU^* \qquad (E\in \E)
\]
for some either unitary or antiunitary operator $U$ on $H$.
Clearly, this result is in an intimate connection with the
fundamental theorem of projective geometry determining the form of
the ortho-order automorphisms of the orthoposet $\mathcal P(H)$ of
all projections on $H$. (In our recent paper \cite{MP} we have
shown that, unlike the fundamental theorem of projective geometry,
the conclusion in Ludwig's result remains valid also in the
two-dimensional case.)

It is an exciting problem to characterize the automorphisms of
algebraic structures of any kind by means of their preservation
properties concerning certain relevant relations, sets,
quantities, etc. which are connected with the underlying
structures. To mention only one such area of investigations, we
refer to the linear preserver problems which, in the last decades,
represent one of the most extensively studied research areas in
matrix theory (see, for example, the survey paper \cite{LiTsing}).
In what follows we present three results of the above kind
concerning the automorphisms of the Hilbert space effect algebra.
Our aim with this paper is to try to draw the attention of the
people working on the foundations of quantum mechanics and dealing
with algebraic structures appearing there to such problems. We
believe that just like in certain parts of pure mathematics, such
investigations can give new insight into the behaviour of the
automorphisms that might help to better understand the underlying
algebraic structures.

We now turn to our results. Let us begin with the following small,
innocent remark.

\begin{remark}
Clearly, Ludwig's theorem describes those bijections of the effect
algebra which preserve the relation $\leq$ and the operation
$\perp$. However, these properties in question can be expressed by
the preservation of one single relation which is the
orthogonality. The effects $E,F$ are said to be orthogonal if
$E\leq I-F$ (or, equivalently, if $E+F\leq I$) (see, for example,
\cite{Gudder2}). Now, our assertion is that a bijective map $\phi:
\E \to \E$ is an ortho-order automorphism of $\E$ if and only if
$\phi$ preserves the orthogonality in both directions. Indeed, the
necessity is obvious. Conversely, suppose that $\phi$ preserves
the orthogonality in both directions. It is easy to see that for
any effects $A,B$ we have $A\leq B$ if and only if for every $C\in
\E$, the orthogonality of $B$ and $C$ implies the orthogonality of
$A$ and $C$. This characterization of the order gives us that
$\phi$ preserves the order in both directions. Next, for any
effect $A\in \E$, the effect $A^\perp$ can easily be characterized
as the supremum of all effects which are orthogonal to $A$. We now
easily get that $\phi$ preserves the operation $\perp$. This
proves our assertion.

We point out that the orthogonality preserving property appears in
the definition of the so-called effect-automorphisms \cite[D
4.2.1]{Ludwig} (in \cite{CVLL2} they were called
{\bf{E}}-automorphisms). These are bijective maps $\phi:\E \to \E$
with the property that for every $E,F\in \E$ we have
\begin{equation}\label{E:E-aut}
E+F\in \E \Longleftrightarrow \phi(E)+\phi(F)\in \E
\end{equation}
and in this case
\[
\phi(E+F)=\phi(E)+\phi(F)
\]
holds. It now follows that the first property \eqref{E:E-aut} in
the definition of effect automorphisms characterizes exactly the
ortho-order automorphisms. Observe that it follows from Ludwig's
theorem if $\dim H\geq 3$ and from \cite{MP} if $\dim H=2$ that
the ortho-order automorphisms are additive, so in those cases
these two kinds of automorphisms are the same. This is trivially
not true if $\dim H=1$.
\end{remark}

We now turn to the nontrivial results of the paper. Beside order
and orthogonality there is another important relation on $\E$.
This is the coexistency (see, for example, \cite[II.2.2.]{BGL} or
\cite[$\S$ 1]{Kraus}). A set of effects is called coexistent if
its members are in the range of an unsharp observable, i.e. a POV
(positive operator valued) measure. In the case of two effects
$E,F$ this is well-known to be equivalent to the following: there
exist effects $A,B,C\in \E$ such that
\[
E=A+C, \quad F=B+C, \text{ and } A+B+C\in \E.
\]

Our first theorem which follows tells us that the preservation of
the two binary relations of order and coexistency characterizes
the ortho-order automorphisms of $\E$.

\begin{theorem}\label{T:coexistency}
Let $H$ be a Hilbert space with $\dim H\geq 3$. Let $\phi: \E \to
\E$ be a bijective map with the properties that
\[
E\leq F \Longleftrightarrow \phi(E)\leq \phi(F)
\]
and
\[
E \text{ and } F \text{ are coexistent } \Longleftrightarrow
\phi(E) \text{ and } \phi(F) \text{ are coexistent }
\]
for every $E,F\in \E$. Then there exists an either unitary or
antiunitary operator $U$ on $H$ such that $\phi$ is of the form
\[
\phi(E)=UEU^* \qquad (E\in \E).
\]
\end{theorem}

We remark that it is easy to see that the preservation of the
order or the preservation of the coexistency alone does not
characterize the automorphisms of $\E$. As for order, see the
remark after the proof of Theorem~\ref{T:coexistency}. As for
coexistency, consider the transformation $\phi:\E \to \E$ defined
as $\phi(0)=I, \phi(I)=0$ and $\phi(E)=E$ otherwise.

If $\varphi$ is a pure state (i.e. a unit vector in $H$), then the
probability of an effect $E\in \E$ in this state is $\la
E\varphi,\varphi\ra$. Our second result asserts that if a
bijective map $\phi$ on $\E$ preserves the order and there are two
pure states $\varphi,\psi\in H$ with respect to which $\phi$
preserves the probability, then $\phi$ is an automorphism of $\E$.
More explicitly, we have the following result.

\begin{theorem}\label{T:probability}
Assume $\dim H\geq 3$. Let $\phi:\E \to \E$ be a bijective map for
which
\[
E\leq F \Longleftrightarrow \phi(E)\leq \phi(F) \qquad (E,F\in \E)
\]
and suppose that there are unit vectors $\varphi, \psi\in H$ such
that
\[
\la \phi(E)\psi,\psi\ra=\la E\varphi,\varphi \ra \qquad (E\in \E).
\]
Then there exists an either unitary or antiunitary operator $U$ on
$H$ such that
\[
\phi(E)=UEU^* \qquad (E\in \E)
\]
and $U\varphi=\psi$.
\end{theorem}

(For further results of the same spirit see the remark after the
proof of Theorem~\ref{T:probability}.) Similarly to the case of
our first theorem, we remark that the preservation of the
probability appearing above alone is not sufficient to
characterize the automorphisms of $\E$. Indeed, choosing pure
states $\varphi=\psi$ and defining $\phi$ by the identity on the
set $\mathcal N$ of all effects $E$ for which $\la
E\varphi,\varphi\ra \neq 0$ and by any permutation on $\E
\setminus \mathcal N$ we can easily get an appropriate example.

The set $\E$ of all effects is clearly a convex set. So, it is
natural to equip it with the operation of convex combinations. The
automorphisms of effect algebras with respect to this operation
which are called mixture automorphisms were studied, for example
in, \cite{Gudder3}. These automorphisms of $\E$ in full generality
were determined in \cite{Molnar}. The result \cite[Corollary
2]{Molnar} says that every mixture automorphism $\phi$ of $\E$ is
either of the form
\[
\phi(E)=UEU^* \qquad (E\in \E)
\]
or of the form
\[
\phi(E)=U(I-E)U^* \qquad (E\in \E),
\]
where $U$ is an either unitary or antiunitary operator on $H$. An
effect $A$ is called a mixture of the effects $B,C$ if $A$ is a
convex combination of $B$ and $C$, that is, if there is a scalar
$\lambda \in [0,1]$ such that $A=\lambda B+(1-\lambda )C$. Our
third result states that the preservation of mixtures
characterizes the mixture automorphisms of $\E$.

\begin{theorem}\label{T:mixture}
Assume $\dim H\geq 2$. Let $\phi:\E \to \E$ be a bijective
function with the property that
\[
A \text{ is a mixture of $B$ and $C$ } \Longleftrightarrow \phi(A)
\text{ is a mixture of $\phi(B)$ and $\phi(C)$}
\]
holds for all $A,B,C\in \E$. Then there exists an either unitary
or antiunitary operator $U$ on $H$ such that either
\[
\phi(A)=UAU^* \qquad (A\in \E)
\]
or
\[
\phi(A)=U(I-A)U^* \qquad (A\in \E).
\]
\end{theorem}

\section{Proofs}

This section is devoted to the proofs of our results. We begin
with some lemmas.

\begin{lemma}\label{L:scalar}
The effect $A\in \E$ is coexistent with every $E\in \E$ if and
only if $A=\lambda I$ holds with some scalar $\lambda \in [0,1]$.
\end{lemma}

\begin{proof}
If $A$ is coexistent with every effect, then it is coexistent with
every projection $P$ on $H$. Since coexistence with a projection
means commutativity with that projection (see, for example,
\cite[p. 120]{Kraus}), it follows that $A$ commutes with every
projection which implies that $A$ commutes with every operator
$B\in B(H)$ ($B(H)$ denotes the algebra of all bounded linear
operators on $H$). It is well-known that this implies that $A$ is
a scalar. Conversely, suppose that $A=\lambda I$ with some
$\lambda \in [0,1]$. If $E\in \E$ is arbitrary, then we can write
\[
\lambda I=\lambda (I-E)+\lambda E, \qquad E=(1-\lambda)E+\lambda
E.
\]
Since $\lambda I+(1-\lambda)E\leq I$, it follows that $\lambda I$
and $E$ are coexistent.
\end{proof}

\begin{lemma}\label{L:rank-1}
Let $E,F\in \E$ be of rank 1. Suppose that the ranges of $E$ and
$F$ are different. Then $E,F$ are coexistent if and only if
$E+F\in \E$.
\end{lemma}

\begin{proof}
Only the necessity requires proof. Suppose that $E,F$ are
coexistent. Then there are effects $A,B,C$ such that $A+B+C$ is an
effect and $E=A+C$, $F=B+C$. As $E,F$ are of rank 1 and $C\leq
E,F$, it follows that the range $\rng C$ of $C$ is included in the
range of $E$ and $F$ which implies that $\rng C=\{ 0\}$, that is,
$C=0$. This shows that $E+F=A+B+2C=A+B+C\in \E$.
\end{proof}

In what follows we need the concept of the strength of effects
along rays (rank-1 projections) defined in \cite{Gudder}. Let
$E\in \E$ and consider an arbitrary rank-1 projection $P$ on $H$.
The strength of $E$ along $P$ is defined by
\[
\lambda(E,P)=\sup\{ \lambda\in [0,1]\, :\, \lambda P\leq E\}.
\]

If $\varphi\in H$ is any unit vector, then let $P_\varphi$ denote
the rank-1 projection which projects onto the linear space
generated by $\varphi$.

In the sequel we shall use the following nice result of Busch and
Gudder \cite[Theorem 3]{Gudder}: for any effect $E\in \E$ and unit
vector $\varphi \in H$ we have
\[
\exists \lambda >0 \, :\, \lambda P_\varphi \leq E\,
\Longleftrightarrow \, \varphi\in \rng E^{1/2}.
\]

\begin{lemma}\label{L:eigenvalue}
Let $E\in \E$ and $0<\lambda<\mu \leq 1$. Suppose that
\[
\lambda I \leq E\leq \mu I.
\]
If $\varphi, \psi \in H$ are unit vectors such that
\[
\lambda(E,P_\varphi)=\lambda, \quad \lambda(E,P_\psi)=\mu,
\]
then $\varphi, \psi$ are eigenvectors of $E$ and the corresponding
eigenvalues are $\lambda,\mu$, respectively.
\end{lemma}

\begin{proof}
It follows from $\lambda I\leq E\leq \mu I$ that for spectrum
$\sigma(E)$ of $E$ we have $\sigma (E)\subset [\lambda,\mu]$. We
assert that $\lambda ,\mu \in \sigma(E)$. Indeed, suppose, for
example, that the effect $E-\lambda I$ is invertible. Then its
square-root is also invertible and the above mentioned result of
Busch and Gudder  (\cite[Theorem 3]{Gudder}) tells us that there
exists a positive number $\epsilon$ for which
\[
\epsilon P_\varphi\leq E-\lambda I.
\]
This implies that
\[
(\epsilon+\lambda )P_\varphi \leq E
\]
which means that the strength of $E$ along $P_\varphi$ is greater
than $\lambda$. But this is a contradiction. So, we have $\lambda
\in \sigma(E)$. One can prove in a similar fashion that $\mu \in
\sigma(E)$. Therefore, the convex hull of $\sigma(E)$ is exactly
$[\lambda, \mu]$. Now, one can follow the proofs of the statements
(a), (b) in \cite[Theorem 5]{Gudder} to verify that
$E\varphi=\lambda \varphi$ and $E\psi=\mu \psi$.
\end{proof}

\begin{remark}
It is easy to see that in the previous lemma $\lambda =0$ can not
be allowed. To show this, pick two different rank-1 projections
$P,Q$ such that $PQ\neq 0$. Let $\varphi \in H$ be a unit vector
such that $Q=P_\varphi$. Clearly, we have $0\leq P\leq I$ and
$\lambda (P,P_\varphi)=0$ but $P \varphi\neq 0\cdot \varphi$.
\end{remark}

\begin{lemma}\label{L:proj}
Let $P,Q$ be different projections on the Hilbert space $H$. Then
$P$ and $Q$ are mutually orthogonal (that is, $PQ=0$) if and only
if every subprojection of $P$ commutes with every subprojection of
$Q$.
\end{lemma}

\begin{proof}
This follows easily form the following observation: the
projections $P,Q$ commute if and only if there are mutually
orthogonal projections $P_0,Q_0,R$ such that $P=P_0+R$ and
$Q=Q_0+R$.
\end{proof}

Now, we are in a position to prove our first theorem.

\begin{proof}[Proof of Theorem~\ref{T:coexistency}]
Any bijection of $\E$ which preserves the order in both directions
also preserves the projections in both directions. This important
observation was made in \cite[Theorem 5.8., p 219]{Ludwig}. By the
order preserving property of $\phi$, one can deduce that the
operator $\phi(P)$ is a rank-1 projection if and only if $P$ is a
rank-1 projection. More generally, one can prove that $\phi(P)$ is
a rank-$n$ projection if and only if $P$ is a rank-$n$ projection.
Indeed, this follows from the following characterization of
rank-$n$ projections: the projection $P$ is of rank-$n$ if and
only if there is a chain $P_1, \ldots ,P_{n-1}$ of $n-1$
projections such that $0\lneq P_1 \lneq P_2 \lneq, \ldots \lneq
P_{n-1} \lneq P$ but there is no such chain of $n$ members.

By Lemma~\ref{L:scalar} we find that there is a bijective strictly
monotone increasing function $f:[0,1] \to [0,1]$ such that
\[
\phi(\lambda I)=f(\lambda)I \qquad (\lambda \in [0,1]).
\]
Let $P$ be a rank-1 projection. Since $\phi(P)$ is also of rank 1,
it follows from $0\leq \phi(\lambda P)\leq \phi(P)$ that
$\phi(\lambda P)$ is a scalar multiple of $\phi(P)$. Therefore, we
have a bijective strictly monotone increasing function
$f_P:[0,1]\to [0,1]$ such that $\phi(\lambda P)=f_P
(\lambda)\phi(P)$ $(\lambda \in [0,1])$. The strength of
$\phi(\lambda I)=f(\lambda )I$ along $\phi(P)$ is obviously
$f(\lambda)$. On the other hand, we have
\[
f_P(\mu) \phi(P)=\phi(\mu P)\leq \phi(\lambda I)=f(\lambda )I
\]
if and only if $\mu \leq \lambda$ which shows that the strength of
$\phi(\lambda I)$ along $\phi(P)$ is $f_P(\lambda)$. Therefore, we
have $f_P=f$. Consequently,
\[
\phi(\lambda P)=f(\lambda) \phi(P) \qquad (\lambda \in [0,1])
\]
holds for every rank-1 projection $P$ on $H$.

Since, as we have already mentioned, the coexistence of
projections is equivalent to the commutativity of the projections
in question, it follows from Lemma~\ref{L:proj} that $\phi$
preserves the orthogonality of projections. So, $\phi$ is a
bijection of the set of all projections on $H$ which preserves the
order and the orthogonality in both directions. It follows from
the fundamental theorem of projective geometry (see the
introduction) that there exists an either unitary or antiunitary
oprerator $U$ on $H$ such that $\phi(P)=UPU^*$ holds for every
projection $P$ on $H$. Considering the transformation $E\mapsto
U^*\phi(E)U$ if necessary, we can clearly assume without any loss
of generality that $\phi(P)=P$ holds for every projection $P$. It
now remains to prove that we have $\phi(E)=E$ for every effect $E$
as well.

Let $P_1, \ldots, P_n$ be pairwise orthogonal rank-1 projections
and $\lambda_1, \ldots, \lambda_n \in [0,1]$. Set $P=P_1 + \ldots
+ P_n$ and $E=\lambda_1 P_1 + \ldots +\lambda_n P_n$. Since
$\phi(E)\leq \phi(P)$, we deduce that $\phi(E)$ acts on the
$n$-dimensional subspace $\rng \phi(P)$ (this means that $\phi(E)$
sends the range of $\phi(P)$ into itself and $\phi(E)$ is zero on
the orthogonal complement of $\rng \phi(P)$). Since each $P_i$
commutes with $E$, it follows from the coexistence preserving
property of $\phi$ that $\phi(P_i)$ commutes with $\phi(E)$. As
the sum of the $\phi(P_i)$'s is $\phi(P)$, we readily obtain that
\begin{equation}\label{E:effect1}
\phi(E)=\mu_1 \phi(P_1)+ \ldots + \mu_n \phi(P_n)
\end{equation}
holds for some scalars $\mu_i\in [0,1]$. Since the strength of $E$
along $P_i$ is $\lambda_i$, by the order preserving property of
$\phi$ we infer that the strength of $\phi(E)$ along $\phi(P_i)$
is $f(\lambda_i)$. On the other hand, it follows from the equality
\eqref{E:effect1} that the strength of $\phi(E)$ along $\phi(P_i)$
is $\mu_i$. Therefore, we have
\begin{equation*}
\phi(\lambda_1 P_1+ \ldots + \lambda_n P_n)=\phi(E)=f(\lambda_1)
\phi(P_1)+ \ldots + f(\lambda_n) \phi(P_n).
\end{equation*}
This gives us that for any finite rank operator $A\in \E$ we have
$\phi(A)=f(A)$, where $f(A)$ denotes the image of $f$ under the
continuous function calculus corresponding to the normal operator
$A$. (Observe that, as $f:[0,1]\to [0,1]$ is a strictly monotone
increasing bijection, it is a continuous function.) It follows
from the spectral theorem and the properties of the spectral
integral that for any effect $A\in \E$ there is a net $(A_\alpha)$
of finite rank effects such that $A_\alpha \leq A$ and $A_\alpha
\to A$ in the strong operator topology. Since the multiplication
is continuous on the bounded subsets of operators with respect to
the strong operator topology, we obtain that $p(A_\alpha) \to
p(A)$ strongly for every polynomial $p$. As, by Weierstrass's
theorem, $f$ can be approximated by polynomials in the uniform
norm, we find that $f(A_\alpha) \to f(A)$ strongly. Since
\[
f(A_\alpha)=\phi(A_\alpha)\leq \phi(A),
\]
we obtain that
\begin{equation}\label{E:effect2}
f(A)\leq \phi(A) \qquad (A\in \E).
\end{equation}
To see the reverse inequality, observe that we have
$\phi^{-1}(\lambda I)=f^{-1}(\lambda)I$. Therefore, considering
$\phi^{-1}$ in the place of $\phi$, we get that
\begin{equation}\label{E:effect3}
f^{-1}(A)\leq \phi^{-1}(A).
\end{equation}
It follows from \eqref{E:effect2} and \eqref{E:effect3} that
\[
A=f^{-1}(f(A))\leq \phi^{-1}(f(A))\leq \phi^{-1}(\phi(A))= A.
\]
This implies that $\phi(A)=f(A)$ holds for every $A\in \E$.

We show that $f(\lambda)+f(1-\lambda)=1$ $(\lambda \in [0,1])$.
Let $0<\lambda <1$. Let $P$ be a rank-1 projection. Pick any
$0<\epsilon <1-\lambda$. Clearly, the spectrum of $\lambda
P+\epsilon P$ is $\{0, \lambda+\epsilon\}$. Let $Q$ be a rank-1
projection such that $\rng Q\cap \rng P=\{ 0\}$. If $\delta$
denotes the largest eigenvalue of the positive operator $\lambda
P+\epsilon Q$, then by Weyl's perturbation theorem (see, for
example, \cite[Corollary III.2.6]{Bhatia}) we have
\[
|(\lambda+\epsilon)-\delta|\leq \| (\lambda P+\epsilon P)-(\lambda
P +\epsilon Q)\| =\epsilon \| P-Q\|.
\]
So, if $Q$ is close enough to $P$, then the largest eigenvalue of
the operator $\lambda P+\epsilon Q$ is close enough to
$\lambda+\epsilon$ and hence it is less then 1 which shows that
$\lambda P +\epsilon Q\in \E$. This implies that the effects
$\lambda P$ and $\epsilon Q$ are coexistent and by the properties
of $\phi$ it follows that the same must hold true for
$f(\lambda)P=f(\lambda)\phi(P)=\phi(\lambda P)$ and
$f(\epsilon)Q=f(\epsilon)\phi(Q)=\phi(\epsilon Q)$. By
Lemma~\ref{L:rank-1} we infer that $f(\lambda)P+f(\epsilon)Q\leq
I$. If we let $Q$ converge to $P$, we get
$f(\lambda)P+f(\epsilon)P\leq I$. This gives us that
$f(\lambda)+f(\epsilon) \leq 1$. Therefore, if $\epsilon$ tends to
$1-\lambda$, we obtain $f(\lambda)+f(1-\lambda) \leq 1$ or,
equivalently,
\begin{equation}\label{E:effect4}
f(1-\lambda) \leq 1-f(\lambda).
\end{equation}
Applying the above argument for $\phi^{-1}$ instead of $\phi$, we
find that
\begin{equation}\label{E:effect5}
f^{-1}(1-\lambda) \leq 1-f^{-1}(\lambda).
\end{equation}
Using the monotonicity of $f^{-1}$ and the inequalities
\eqref{E:effect4}, \eqref{E:effect5}, we have
\[
\lambda=f^{-1}(f(\lambda))\leq f^{-1}(1-f(1-\lambda))\leq
1-f^{-1}(f(1-\lambda))=\lambda.
\]
Therefore, we deduce $f^{-1}(f(\lambda))= f^{-1}(1-f(1-\lambda))$
which implies that
\[
f(\lambda)+f(1-\lambda)=1
\]
for every $0<\lambda<1$. It is trivial that the equality is valid
for $\lambda=0,1$ as well. This implies that
\[
f(I-A)=I-f(A)
\]
holds for every effect $A\in \E$ which yields that $\phi$
satisfies
\[
\phi(I-A)=I-\phi(A) \qquad (A\in\E).
\]
Therefore, $\phi$ is an ortho-order automorphism of $\E$. Applying
Ludwig's theorem on the form of those automorphisms we have
$\phi(E)=E$ $(E\in \E)$. This completes the proof of the theorem.
\end{proof}

\begin{remark}
A careful examination of the proof of Theorem~\ref{T:coexistency}
shows that if $\phi: \E \to \E$ is a bijective map which preserves
the order and the commutativity in both directions, then there
exists an either unitary or antiunitary operator $U$ on $H$ and a
strictly monotone increasing bijection $f:[0,1]\to [0,1]$ such
that $\phi$ is of the form
\begin{equation}\label{E:effect6}
\phi(A)=Uf(A)U^* \qquad (A\in \E).
\end{equation}
Clearly, it follows from the order preserving property of $\phi$
that $f$ as well as $f^{-1}$ are operator monotone on $[0,1]$.
(Recall that a continuous real function $g$ on an interval is
called operator monotone if for arbitrary selfadjoint operators
$A, B$ with spectrum in the domain of $g$, the relation $A\leq B$
implies $g(A)\leq g(B)$.)

It is easy to see that, conversely, if $U$ and $f$ are such as
above, then the formula \eqref{E:effect6} defines a bijection of
$\E$ which preserves the order and commutativity in both
directions. (This last property follows from the fact that if
$A,B$ are commuting, then the same holds for their polynomials.
Finally, as every continuous function on a compact subset of the
real line can be uniformly approximated by polynomials, we obtain
the commutativity of any continuous function of $A$ and $B$.) Now,
the question is that whether there do exist nontrivial continuous
functions $f:[0,1]\to [0,1]$ with the property that $f,f^{-1}$ are
operator monotone. The answer to this question is affirmative.
Indeed, consider, for example, the function
$f(\lambda)=(2\lambda)/(1+\lambda)$ $(\lambda \in [0,1])$.

It is a remarkable fact from the mathematical point of view that
the preservation of the order and the coexistency together do
characterize the automorphisms of $\E$ while the preservation of
the order and the commutativity (which is the closest widely used
property in pure mathematics to coexistency) do not.
\end{remark}

We continue with the proof of our second theorem. We shall need
the following observation.

\begin{lemma}\label{L:scalar2}
Let $E\in \E$ and let $D$ be a dense subset of the set of all unit
vectors in $H$. Pick $0<\lambda \leq 1$. If $\lambda
(E,P_\varphi)=\lambda$ for every $\varphi \in D$, then $E=\lambda
I$.
\end{lemma}

\begin{proof}
As $\lambda P_\varphi \leq E$ for every $\varphi \in D$ and $D$ is
dense in the set of all unit vectors, it follows that $\lambda
P_\varphi\leq E$ holds for every unit vector $\varphi$ in $H$.
According to the result \cite[Theorem 3]{Gudder} we deduce that
the square-root of $E$ is surjective which gives us that $E^{1/2},
E$ are invertible. In \cite[Theorem 4]{Gudder}, Busch and Gudder
gave an explicit formula for the strength of an arbitrary effect
along an arbitrary ray. It follows form that result that
$\|E^{-1/2}\varphi\|^{-2}=\lambda$ for every $\varphi\in D$ which
implies that we have the same equality for every unit vector
$\varphi\in H$. This gives us that
\[
\la E^{-1}\varphi,\varphi\ra =\la \frac{1}{\lambda} \varphi,
\varphi \ra
\]
for every $\varphi \in H$. Hence, we have $E=\lambda I$.
\end{proof}

We now can prove Theorem~\ref{T:probability}.

\begin{proof}[Proof of Theorem~\ref{T:probability}]
Just as in the proof of Theorem~\ref{T:coexistency}, we obtain
that $\phi$ preserves the projections in both directions as well
as their rank and that for every rank-1 projection $P$ there is a
strictly monotone increasing bijection $f_P:[0,1]\to [0,1]$ such
that $\phi(\lambda P)=f_P(\lambda) \phi(P)$ $(\lambda \in [0,1])$.

Let $\varphi,\psi$ be as in the theorem. If the range of $P$ is
not orthogonal to $\varphi$, then we compute
\[
f_P(\lambda)\la \phi(P)\psi,\psi\ra =\la \phi(\lambda
P)\psi,\psi\ra= \la \lambda P\varphi,\varphi\ra= \lambda \la
P\varphi,\varphi \ra.
\]
This gives us that $f_P(\lambda)=c\lambda$ $(\lambda \in [0,1])$
for some constant $c$. Since $f_P$ is a bijection of $[0,1]$ onto
itself, it follows that $c=1$. So, in the present case we have
\[
\phi(\lambda P)=\lambda \phi(P) \qquad (\lambda \in [0,1]).
\]

Since
\[
\| \phi(P_\varphi)\psi\|^2= \la \phi(P_\varphi)\psi,\psi\ra=\la
P_\varphi \varphi,\varphi\ra=1,
\]
we deduce that $\| \phi(P_\varphi)\psi\|=1=\| \psi\|$ which
implies that $\psi$ is in the range of $\phi(P_\varphi)$. Hence we
have $\phi(P_\varphi)=P_\psi$.

Suppose now that $PP_\varphi=0$. As
\[
\la \phi(P)\psi,\psi\ra=\la P\varphi,\varphi\ra=0,
\]
we have $\phi(P)\psi=0$ implying that $\phi(P)P_\psi=0$.
Therefore, $\phi(P)$ is orthogonal to $P_\psi=\phi(P_\varphi)$. As
$\phi^{-1}$ has similar properties as $\phi$, we obtain that a
rank-1 projection $P$ is orthogonal to $P_\varphi$ if and only if
$\phi(P)$ is orthogonal to $\phi(P_\varphi)$.

Let $P$ be a rank-1 projection orthogonal to $P_\varphi$. Let
$0<\lambda \leq 1$ be arbitrary but fixed and consider the
operator $A=\phi(\lambda(P+P_\varphi))$. Set
$Q=\phi(P+P_\varphi)$. Since $A\leq Q$, it follows that $A$ acts
on the range of the rank-2 projection $Q$. Let $Q_0$ be any rank-1
subprojection of $Q$ which is not orthogonal to $\phi(P_\varphi)$.
Then $\phi^{-1}(Q_0)$ is a rank-1 subprojection of $P+P_\varphi$
which is not orthogonal to $P_\varphi$. Clearly, the strength of
$\lambda (P+P_\varphi)$ along $\phi^{-1}(Q_0)$ is $\lambda$. It
follows from the second section of the present proof that the
strength of $A$ along $Q_0$ is also $\lambda$. Since $Q_0$ runs
through the set of all rank-1 subprojections of $Q$ which are not
orthogonal to $\phi(P_\varphi)$, Lemma~\ref{L:scalar2} applies to
obtain $A=\lambda Q$. Therefore, we have
\[
f_P(\lambda)\phi(P)=\phi(\lambda P)\leq
\phi(\lambda(P+P_\varphi))=A=\lambda Q.
\]
This gives us that $f_P(\lambda)\leq \lambda$ $(\lambda \in
[0,1])$. (Observe that in fact we have the above inequality only
for positive $\lambda$'s but $f_P(0)=0$ is trivial because of the
definition of $f_P$.) Applying the above argument in relation with
$\phi^{-1}$ in the place of $\phi$, we find that $f_P^{-1}(\lambda
)\leq \lambda$ $(\lambda\in [0,1])$. Since, due to the order
preserving property of $\phi$, $f_P$ is monotone increasing, it
follows that $f_P(\lambda)=\lambda$.

To sum up what we have already proved, we have
\begin{equation}\label{E:effect9}
\phi(\lambda P)=\lambda \phi(P)
\end{equation}
for every $\lambda \in [0,1]$ and rank-1 projection $P$ on $H$ no
matter $P$ is orthogonal to $P_\varphi$ or not. By
\eqref{E:effect9} the strength of $\phi(\lambda I)$ along every
rank-1 projection is $\lambda$ and hence, by
Lemma~\ref{L:scalar2}, we have $\phi(\lambda I)=\lambda I$.

Now we are in a position to prove that $\phi$ preserves the
orthogonality between projections. Let $P,Q$ be rank-1 projections
with $PQ=0$. Choose a projection $R$ such that $P\leq R$ and
$Q\leq I-R$. Let $0<\lambda< \mu \leq 1$. Set $E=\lambda R+\mu
(I-R)$. We have
\[
\lambda I=\phi(\lambda I)\leq \phi(E)\leq \phi(\mu I) \leq \mu I.
\]
Since the strength of $E$ along $P$ is $\lambda$, it follows from
\eqref{E:effect9} that the strength of $\phi(E)$ along $\phi(P)$
is also $\lambda$. Similarly, we obtain that the strength of
$\phi(E)$ along $\phi(Q)$ is $\mu$. Lemma~\ref{L:eigenvalue} shows
that $\rng \phi(P),\rng \phi(Q)$ are eigensubspaces of $\phi(E)$
and the corresponding eigenvalues are $\lambda, \mu$,
respectively. Since the eigensubspaces of a self-adjoint operator
corresponding to different eigenvalues are mutually orthogonal, it
follows that $\phi(P)\phi(Q)=0$. Since $\phi$ is order-preserving
and every projection is the supremum of all rank-1 projections
which are included in it, this implies that $\phi$ preserves the
orthogonality between arbitrary projections. So, for any
projections $P,Q$ on $H$ we have $PQ=0$ if and only if
$\phi(P)\phi(Q)=0$.

It follows that $\phi$ is a bijection of the set of all
projections on $H$ which preserves the order and the orthogonality
in both directions. By the fundamental theorem of projective
geometry there exists an either unitary or antiunitary oprerator
$U$ on $H$ such that
\[
\phi(P)=UPU^*
\]
holds for every projection $P$ on $H$. By \eqref{E:effect9} it
follows that
\begin{equation}\label{E:weak}
\phi(\lambda P)=U(\lambda P)U^*
\end{equation}
for every $\lambda \in [0,1]$ and rank-1 projection $P$. The
operators of the form $\lambda P$ are the weak atoms in $\E$
\cite[Lemma 2]{Gudder}. The statement \cite[Corollary 3]{Gudder}
says that every effect $E$ is the supremum of all weak atoms which
are less than or equal to $E$. It now follows form \eqref{E:weak}
that
\[
\phi(A)=UAU^* \qquad (A\in \E).
\]
If $U$ above is unitary, then we have
\[
\la AU^*\psi,U^*\psi\ra =\la \phi(A)\psi,\psi\ra=\la
A\varphi,\varphi \ra
\]
for every $A\in \E$. It is easy to see that this implies $U^*\psi
=\epsilon \varphi$ for some complex number $\epsilon$ of modulus
1. This yields $\epsilon U\varphi=\psi$. Replacing $U$ by
$\epsilon U$, we obtain the last assertion of our theorem. If $U$
is antiunitary, then one can argue is a similar way. This
completes the proof of the theorem.
\end{proof}

\begin{remark}
We note that the same conclusion as in Theorem~\ref{T:probability}
holds true if $\phi:\E \to \E$ is a bijective map which preserves
the order in both directions as well as the spectrum (or the
numerical range, or, more generally, the numerical radius = the
spectral radius = the norm of effects). We omit the proofs since
those results have no real physical content (with the only
possible exception of the numerical range which concept can be
interpreted as the set of all probabilities of an effect
corresponding to pure states).
\end{remark}

We now turn to the proof of our last theorem. Just as before, we
need some auxiliary results.

\begin{lemma}\label{L:bounded}
Let $f$ be a linear functional on the real linear space $B_s(H)$
of all bounded self-adjoint operators on $H$. If $f$ is bounded
from below on the operator interval $[0,I]$, then $f$ is a bounded
linear functional.
\end{lemma}

\begin{proof}
Let $K$ be a real number such that $K\leq f(A)$ $(A\in [0,I])$. We
show that $f$ is bounded also from above on $[0,I]$. Suppose on
the contrary that for every $n\in \N$ there exists an operator
$A_n \in [0,I]$ such that $f(A_n)\geq 2^n$. Let
$B=\sum_{n=1}^\infty A_n/2^n \in [0,I]$ and $B_n =\sum_{k=1}^n
A_k/2^k\in [0,I]$. Clearly, we have $B-B_n\in [0,I]$ and hence
$K\leq f(B-B_n)$ which implies that
\[
K+n\leq K+f(B_n)\leq f(B).
\]
Since this holds for every $n\in \N$ we arrive at a contradiction.
This gives us that $f$ is bounded on $[0,I]$. Since every
self-adjoint operator of norm not greater than 1 is the difference
of two elements of $[0,I]$, we obtain the boundedness of $f$.
\end{proof}

\begin{lemma}\label{L:zero}
Assume $\dim H\geq 2$. Let $f$ be a bounded linear functional on
$B_s(H)$. Suppose that $f(P)=c$ for every nonzero projection $P$
where $c$ is a fixed scalar. Then we have $f=0$.
\end{lemma}

\begin{proof}
Clearly, every projection $P\neq I$ is the difference of two
nonzero projections. We thus obtain that $f$ is 0 on the set of
all such projections. Since $I$ is the sum of two projections
different from $I$, we obtain that $f$ vanishes on the whole set
of projections. By spectral theorem, the linear span of all
projections is norm-dense in $B_s(H)$ and hence we obtain that
$f=0$.
\end{proof}

\begin{proof}[Proof of Theorem~\ref{T:mixture}]
There is a beautiful result due to P\'ales \cite{Pales} on segment
preserving maps between general convex sets in linear spaces. Its
assertion can be translated to our situation in the following way:
if $K$ is a noncollinear convex set in a real linear space $X$ and
$\phi:K\to K$ is a bijective function with the property that $x$
is a mixture (i.e., a convex combination) of $y,z$ if and only if
$\phi(x)$ is a mixture of $\phi(y),\phi(z)$ $(x,y,z\in K)$, then
$\phi$ can be written in the form
\begin{equation}\label{E:zsolt}
\phi(x)=\frac{\psi(x)+ b}{f(x)+c} \qquad (x\in K),
\end{equation}
where $\psi:X\to X$ is a linear transformation, $b\in X$ is fixed,
$f:X\to \mathbb R$ is a linear functional, $c\in \R$ is fixed, and
the denominator in \eqref{E:zsolt} is everywhere positive on $K$.

Adapting this result for $\E$, we have a linear transformation
$\psi$ on $B_s(H)$, an operator $B\in B_s(H)$, a linear functional
$f:B_s(H)\to \R$ and a constant $c\in \R$ such that $f+c$ is
positive on $\E$ and
\begin{equation}\label{E:mix2}
\phi(A)=\frac{\psi(A)+B}{f(A)+c} \qquad (A\in \E).
\end{equation}
By Lemma~\ref{L:bounded}, $f$ is a bounded linear functional on
$B_s(H)$.

Since $0\leq \phi(A)\leq I$ for every $A\in \E$, it follows that
$0\leq \psi(A)+B\leq (f(A)+c)I$ $(A\in \E)$. If $M>0$ denotes an
upper bound of the values of $f+c$ on $\E$, then we have
\[
-B \leq \psi(A)\leq MI -B \qquad (A\in \E).
\]
This implies that
\[
-\|B\| I\leq \psi(A)\leq (M + \|B\|)I \qquad (A\in \E)
\]
which yields that the numerical range of the operator $\psi(A)$
$(A\in \E)$ is contained in the interval $[-\|B\|, M+\| B\|]$.
Since the numerical radius and the norm of a selfadjoint operator
coincide, we obtain that
\[
\| \psi(A)\| \leq M+\| B\| \qquad (A\in \E),
\]
that is, $\psi$ is bounded on $\E$. Just as in the proof of
Lemma~\ref{L:bounded} this implies the boundedness of the linear
transformation $\psi$. By the continuity of $\psi$ and $f$ we
obtain that $\phi$ is norm-continuous and, as $\phi^{-1}$ has the
same properties as $\phi$, we deduce that $\phi^{-1}$ is also
continuous, that is, $\phi$ is a homeomorphism of $\E$.

By the preserving property of $\phi$ it follows that $\phi$
preserves the extreme points of the convex set $\E$ in both
directions. It is well-known that the extreme points of $\E$ are
exactly the projections, hence we obtain that $\phi$ preserves the
projections in both directions. The two trivial projections $0,I$
are distinguished in the set $\mathcal P(H)$ of all projections by
the following property: $0,I$ are the only projections which
cannot be connected to a different projection via a continuous
curve inside the set of all projections (see the proof of
\cite[Theorem 1]{Molnar}). By the known properties of $\phi$ we
infer that $\phi$ permutes the projections $0,I$, that is, it maps
0 either to 0 or to $I$. Considering the transformation $A\mapsto
\phi(I-A)$ if necessary, we can assume that $\phi$ sends 0 to 0
and $I$ to $I$.

Since $\psi$, being a linear transformation, maps 0 to 0, it now
follows from \eqref{E:mix2} that the operator $B$ is 0. So, we
have
\begin{equation}\label{E:mix3}
\phi(A)=\frac{\psi(A)}{f(A)+c} \qquad (A\in \E).
\end{equation}
Since $\phi$ sends projections to projections, it follows from
\eqref{E:mix3} that the linear transformation $\psi$ sends every
projection to a scalar multiple of a projection which scalar
might, of course, depend on the projection in question. We show
that in fact there is no such dependence. Let $P$ be a nontrivial
projection. We have nonzero projections $P',Q'$ and nonzero
scalars $\lambda,\mu,\nu$ such that
\[
\psi(P)=\lambda P', \quad \psi(I-P)=\mu Q', \quad \psi(I)=\nu I.
\]
Observe that $\nu$ does not depend on $P$. By the additivity of
$\psi$ we obtain
\[
\nu I=\lambda P'+\mu Q'
\]
and this implies
\begin{equation}\label{E:mix4}
\lambda P'= \nu I-\mu Q'=\nu (I-Q') +(\nu-\mu)Q'.
\end{equation}
Clearly, $Q'\neq I$. From the equation \eqref{E:mix4} we then
easily infer that $\lambda=\nu$. This shows that $(1/\nu)\psi(P)$
is a projection for every nonzero projection $P$.

Since $\phi$ sends projections to projections, we easily obtain
from \eqref{E:mix3} that $\nu/(f(P)+c)=1$ for any nonzero
projection $P$. Since $\nu, c$ are constants,  we infer from
Lemma~\ref{L:zero} that $f=0$. We have $c=\nu$. Therefore, by
\eqref{E:mix3} we conclude that $\phi=(1/\nu) \psi$ holds on $\E$
which shows that $\phi$ extends to a linear transformation on
$B_s(H)$. Therefore, $\phi$ is a mixture automorphism of $\E$ and
\cite[Corollary 2]{Molnar} applies to complete the proof.
\end{proof}

\begin{remark}
We mention that Theorem~\ref{T:mixture} can be generalized for the
case of von Neumann algebras. Namely, one can easily modify the
proofs of Lemma~\ref{L:bounded} and Lemma~\ref{L:zero} as well as
Theorem~\ref{T:mixture} (one should also consult the proof of
\cite[Theorem 1]{Molnar}) to obtain the following statement: If
$\mathcal A\neq \mathbb C I$ is a von Neumann factor on the
Hilbert space $H$ and $\phi$ is a bijection of the effect algebra
of $\mathcal A$ (which is the convex set $[0,I]\cap \mathcal A$)
that preserves mixtures in both directions, then $\phi$ is a
mixture automorphism.

As for Theorem~\ref{T:coexistency} and
Theorem~\ref{T:probability}, the presented proofs heavily depend
on the fact that every projection on the underlying Hilbert space
belongs to the effect algebra. In our opinion, it would be a nice
achievement if one could get any extension of those theorems for
the case of effect algebras of von Neumann algebras.
\end{remark}

\bibliographystyle{amsplain}

\end{document}